\documentclass[12pt]{amsart}
\usepackage{amsmath, amssymb, amsfonts}

\newtheorem {Theorem}{Theorem}[section]

\newtheorem {Proposition}[Theorem]{Proposition}

\newtheorem {theorem}{Theorem}[section]
\newtheorem {lemma}[theorem]{Lemma}

\newtheorem {remark}[theorem]{Remark}

\newcommand{\di}{\mbox{div}}

\newcommand{\R}{\mathbb{R}}
\newcommand{\p}{\partial}
\newcommand{\ba}{\begin{array} }
\newcommand{\ea}{\end{array} }

\def\R{{\mathbb{R}}}

\newcommand{\be}{\begin }

\numberwithin{equation}{section}

\def\ba{\begin{aligned}}
\def\ea{\end{aligned}}
\def\be{\begin{equation}}
\def\ee{\end{equation}}
\def\tF{F}
\def\leqs{\lesssim}

\def\t#1{\tilde{#1}}
\def\eps{\varepsilon}
\def\ts{\tau}

\def\r{\rho}

\def\div{\text{div}}

\def\O{\Omega}

\def\f{\frac}
\def\p{\partial}
\def\a{\alpha}
\def\g{\nabla}
\def\ld{\lambda}
\def\lap{\triangle}

\def\B{\mathbf{B}}
\def\u{\mathbf{u}}
\def\v{\delta}

\def\O{\Omega}

\allowdisplaybreaks

\begin{document}
	\vspace{1cm}
\title[Global strong solutions to 2d MHD]{Global strong solutions to the radially symmetric compressible MHD equations in 2D solid balls with arbitrary large initial data}
\author{Xiangdi Huang, Wei YAN$^*$}
\thanks{* Corresponding author.}

\address{Xiangdi HUANG\hfill\break\indent
	Institute of Mathematics,
\hfill\break\indent
	Academy of Mathematics and Systems Sciences,\hfill\break\indent
	Chinese Academy of Sciences, Beijing 100190, China}
\email{xdhuang@amss.ac.cn}
\address{Wei YAN\hfill\break\indent
	School of Mathematics,\hfill\break\indent
	Jilin Unversity, Changchun 130012, China}
\email{wyanmath@jlu.edu.cn}
\maketitle

\begin{abstract}

In this paper, we prove the global existence of strong solutions to the two-dimensional compressible MHD equations with density dependent viscosity coefficients (known as Kazhikhov-Vaigant model) on 2D solid balls with arbitrary large initial smooth data where shear viscosity $\mu$ being constant and the bulk viscosity $\lambda$ be a polynomial of density up to power $\beta$.

The global existence of the radially symmetric strong solutions was established under Dirichlet boundary conditions for $\beta>1$.  Moreover, as long as $\beta\in (\max\{1,\frac{\gamma+2}{4}\},\gamma]$, the density is shown to be uniformly bounded with respect to time. This generalizes the previous result of \cite{2022Li,2016Huang,2022Huang} of the compressible Navier-Stokes equations on 2D bounded domains where they require $\beta>4/3$ and also improves the result of \cite{chen2022,2015Mei-1,2015Mei-2}  of compressible MHD equations  on 2D solid balls.

\end{abstract}

\vspace*{10pt}

\section{Introduction and main results}

This paper is concerned with the compressible MHD system in a bounded domain $\O\subset \R^2$: 
\begin{equation} \label{eq1}
	\left\{\begin{array}{lr} 
		\rho_t +\di(\rho \u)=0,\\
		(\rho \u)_t+\di(\rho \u \otimes \u)+ \nabla P=\mu\Delta\u+\nabla ((\mu+\lambda)\di \u)-\g(\f 12 |\B|^2) + \B\cdot\g \B,\\
		\B_t + u\cdot\g \B - \B\cdot\g u + \B\div u = \lap \B,\\
		\div \B=0.
	\end{array}\right.
\end{equation}
with the initial data
\begin{equation}\label{initial data}
\rho(x,0)=\rho_{0}(x),\ \u(x,0)=\u_{0}(x),\ \B(x,0)=\B_{0}(x)~~~~\text{ for } x\in \O,
\end{equation}
and Dirichlet boundary conditions
\begin{equation}\label{BD}
\u(x,t)=0,\quad \B(x,t)=0~~~~~on ~~\p\O.
\end{equation}
Here $x\in\O$ is the spatial coordinate, $t\geq 0$ is time. $\rho \geq 0$, $\u$ and $P$ denote the fluid density, velocity and pressure, respectively.

Without lossing of generality, we make the following assumptions.

$\mu$ and $\lambda$ are the viscosity coefficients satisfying the following restrictions:
\begin{equation}
\mu= {\mbox const> 0},\; \lambda(\rho)=b\rho^{\beta},
\end{equation}
with $b>0, ~\beta>1.$  We consider the isentropic gases for which the equation of state is given by 
\begin{equation*}\label{pressure}
P(\rho)=A\rho^{\gamma},
\end{equation*}
where $A$ is a positive constant, $\gamma>1$ is the adiabatic exponent. 

We are interested in the strong solution of the compressible MHD equations which has a long history.

For the classical model with constant viscosity coefficients, the local existence of strong solutions were established by  Vol’pert–Hudjaev\cite{1972Vol} with positive density and Fan-Yu\cite{2009Fan} with nonnegative density. The global classical solution was established by Kawashima\cite{1983Ka}  as long as the initial data is a small perturbation of a non-vacuum constant. By following the idea and proof of Huang-Li-Xin\cite{huang2012global} for the compressible Navier-Stokes equations, Li-Xu-Zhang\cite{2013lihl}  proved the global existence of classical solutions to the Cauchy problem of the MHD equations in the presence of vacuum.

It is also natural to consider the models with variable viscosity coefficients. The difficulties inherit from the classical compressible Navier-Stokes equations. Bresch and Desjardins \cite{BD1} proposed a new entropy inequality (BD-entropy) under an additional constraint on the viscosity coefficients, which played an important role in proving the existence of weak solutions. Such an entropy pair requires a very special relationship between the shear viscosity $\mu$ and the bulk viscosity $\lambda$. 

Another breakthrough is due to Kazhikhov and Vaigant\cite{1995Kazhikhov} where they proposed the assumption that  $\mu$ being a positive constant and $\lambda=\rho^\beta$.  Then they proved the global existence of strong solution on two-dimensional periodic domains for the compressible Navier-Stokes equations with arbitrary large initial data under the condition $\beta>3$ when the density is absence of vacuum. Recently, by employing a new structure and regularity criterion, the condition on $\beta$ was relaxed to $\beta>4/3$ , Huang and Li  \cite{2016Huang} also established the global existence of classical solutions to the periodic domain and  the whole space \cite{2022Huang} with vacuum at far field in presence of vacuum. One can also refer to Jiu-Wang-Xin \cite{JWX} for the Cauchy problem with positive density at far field. Very recently, by using Riemann mapping theorem and the pull-back Green’s function method applied to the effective viscous flux, Fan-Li-Li\cite{2022Li} showed the same result on bounded domains under Navier-slip boundary condition for  $\beta>4/3$.

Let us turn back to the compressible MHD system with density-dependent viscosity. Lv and Huang [31] obtained the local strong solutions of the Cauchy problem of the 2D compressible MHD equations with a vacuum as far field density. Under the condition $\beta>\frac{4}{3}$, Mei \cite{2015Mei-1,2015Mei-2} proved the global existence  of classical solutions to the 2D compressible MHD equations with large initial data and vacuum on the torus $T^2$ and the whole space $R^2$. Motivated by Fan-Li-Li, Chen-Huang-Shi\cite{2022Shi} successfully extend the previous result to the compressible MHD equations under Navier-slip boundary condition with $\beta>\frac{4}{3}$.

In this manuscript, we would like to improve the result of Chen-Huang-Shi for 2D solid balls under Dirichlet boundary conditions. This requires more delicate analysis and new observations to dealing with the boundary term arising from commutator's estimates which was first found by Huang-Su-Yan-Yu\cite{2023Huang} recently. 

In this paper, we consider the global existence and large time behavior of spherically symmetric strong solutions to the initial boundary value problem for the compressible MHD system \eqref{eq1}. 

Let $\O$ be a ball of radius $R$ centered at the origin in $\R^2$, and $A$ and $b$  be unity for convenience. We are concerned with spherically symmetric strong solutions to system $\eqref{eq1}$ with initial data $\eqref{initial data}$ and Dirichlet boundary condition $\eqref{BD}$. 

To this end, we denote 
\begin{equation}\label{sym}
	|x|=r,~~\rho(x,t)=\rho(r,t),~~\u(x,t)=u(r,t)\dfrac{x}{r},~~\B(x,t)=B(r,t)\dfrac{x^{\perp}}{r}.
\end{equation}
The equations  $\eqref{eq1}-\eqref{BD}$ is then transformed to
\begin{equation} \label{eq2}
	\left\{\begin{array}{lr} 
		\p_t\rho_t +\p_r(\rho u)+\dfrac{\rho u}{r}=0,\\
		\p_t(\rho u)+\p_r(\rho u^2)+\dfrac{\rho u^2}{r}+ \p_r P=\p_r[(2\mu+\rho^\beta)(\p_r u+\dfrac{u}{r})] -\p_r(\f 12 B^2) -\f{B^2}{r},\\
		\p_rB_t + u\p_r B + B\p_r u = \p_{rr}B + \f{1}{r}\p_r B -\f{1}{r^2}B,\\
	\end{array}\right.
\end{equation}
with the initial data
\begin{equation}\label{initial data2}
	\rho(r,0)=\rho_{0}(r), u(r,0)=u_{0}(r), B(r,0)=B_0(r),~~~~0\leq r<R,
\end{equation}
and Dirichlet boundary condition
\begin{equation}\label{BD2}
	u(0,t)=u(R,t)=0, B(0,t)=B(R,t)=0,~~~~~t>0.
\end{equation}

The first result concerning the global existence of strong solutions is described as follows:
\begin{theorem}\label{them1}
	Suppose that 
	\be\label{cdt}
	 \beta>1,~ \gamma>1,
	 \ee 
	 and the initial data $(\rho_0,\u_0)$ satisfy 
	\be \label{HD0}
	0\leq \rho_0\in W^{1,q}(\O),\quad \u_0, \B_0\in H^1(\O),\quad q>2,
	\ee
	and satisfy the following compatibility condition
	~\be\label{HD}
	-\mu\Delta \u_0 -\nabla ((\mu+\lambda(\rho_0)) \ div \u_0) + \nabla P(\rho_0) +\f 12\g |\B_0|^2 - \B_0\cdot \g \B_0 = \sqrt{\rho_{0}} g,
	\ee
	where ~$g\in L^2(\Omega).$ Then the problem \eqref{eq1}--\eqref{BD} has a unique strong solution satisfying 
	\be
	\rho\in C([0,T); W^{1,q}),\quad \u\in C([0,T);  H^1),
	\ee  and ~\be\left\{
	\ba
	&0\leq \rho(x,t)\leq C(T),\\
	&\|\sqrt{\rho}\dot \u\|_{L^\infty L^2} + \|\g\dot \u\|_{L^2L^2}\leq C(T),\\
	&\|\g \B\|_{L^\infty L^2} +\|\g^2\B\|_{L^2 L^2}\leq C(T).
	\ea\right.\ee
\end{theorem}

A few remarks are in order:
\begin{remark}
The most problematic term lies in a boundary integral arising from some commutator's estimates, analogous to our previous result \cite{2023Huang} for the compressible Navier-Stokes equations. So the main idea is to follow the key proof of \cite{2023Huang}. Luckily, the magnetic field also works well due to its dissipative structure. This will be the key assumption to ensure the global existence. Indeed, in the forth coming paper \cite{2023HXY}, we will show the strong solution to the system \ref{eq1} will develop finite time singularity as long as the magnetic field is absence of  dissipation.
\end{remark}

Next, in the following Theorem, we will show the large time behavior and uniform boundness of density with respect to time when $\beta\in(\max\{1,\frac{\gamma+2}{4}\},\gamma]$.

\begin{theorem}\label{them2}
Assume that
\be
\gamma>1,~\max\{1,\frac{\gamma+2}{4}\}<\beta\le\gamma.
\ee
Then there is a positive constant C depending only on $\mu,\beta,\gamma,\|\rho_0\|_{L^\infty}$ and $\|\u_0\|_{H^1}$ such that
\be
\sup_{0\le t<\infty}\|\rho\|_{L^\infty}\le C
\ee
and the following asymptotic behaviors hold
\be
\lim_{t\rightarrow\infty}(\|\rho-\rho_s\|_{L^p}+\|\nabla \u\|_{L^p} +\|\B\|_\infty + \|\g\B\|_2)=0
\ee
for any $p\in[1,\infty)$, where
\be
\rho_s=\frac{1}{|\O|}\int_{\Omega}\rho_0dx.
\ee
\end{theorem}

\begin{remark}
The result is parallel to \cite{2023Huang}.
\end{remark}

\section{Preliminaries}
In order to prove Theorem \ref{them1} and \ref{them2}, we will frequently employ the following lemmas.
\begin{lemma}\cite{Kazhikhov}\label{ME} For any~$ q >2$, there is a constant $C_q>0,$ such that for any function $\u(x)$, if~$\u\big|_{\p\O}=0$ or~$\int_\O \u dx = 0$, then
	\be
	\| \u\|_{q} \leq C_q \| \nabla \u \| _{\frac{2q}{q+2}}.
	\ee
\end{lemma}

\begin{lemma}\label{lem2}
	Assume that~$\u$ and~$\B$ satisfy \eqref{sym} and \eqref{BD2}, then it holds that
	~\be
	\| u \|_{\infty} \leqs \| \div \u \|_{2}=\| \nabla \u \| _{2}.
	\ee
	and
	~\be
	\| \B \|_{\infty} \leqs\| \nabla \B \| _{2}.
	\ee
\end{lemma}
\proof
For any $r\in (0,R),$ directly calculations lead to
\[\ba
|u(r)|^2 &= 2\int_0^r u\p_r u dr \leq 2\int_0^R|\p_r u||\f{u}{r}| r dr \\
&\leq \int_0^R|\p_r u|^2 r dr + \int_0^R|\f{u}{r}|^2 r dr = \int_0^R \Big(\p_r u + \f{u}{r}\Big)^2 rdr
=\|\div \u\|_2^2,
\ea\]
which implies that~$ \| u \| _{\infty} \leq \| \div \u \|_{2}$.
The inequality for $\B$ can be proved similarly.
\endproof

The following Sobolev inequality will be used frequently.
\begin{lem}\label{lem3} (See \cite{Kazhikhov}) If~$\u|_{\p\O}=0$ or~$\int_\O \u dx = 0$, then for any ~$q>2$, we have
	~\be\label{Sob1}
	\|\u\|_q\leq C\|\u\|_2^{\f{2}{q}} \|\g \u\|_2^{1-\f{2}{q}},\quad \|\u\|_q \leq C(\|\u\|_2 + \|\g \u\|_2),
	\ee
	and
	~\be\label{Sob}
	\|\u\|_\infty\leq C\|\u\|_2^{\f{q-2}{2q-2}} \|\g \u\|_q^{\f{q}{2q-2}}.
	\ee
\end{lem}

The following Beale--Kato--Majda type inequality can be found in \cite{BKM}.
\begin{lem}(See \cite{BKM})
	For any~$2< q < \infty$, there is a constant~$C(q)$, such that the following estimate holds for all ~$\nabla \u \in W^{1,q}(\Omega)$,
	\be\label{CL}
	\| \nabla \u \|_{\infty} \leq C(\| \div \u \|_{\infty}+\| curl \u \|_{\infty} )\log(e + \| \nabla ^{2} \u \|_{L^q}) + C\| \nabla \u \|_{L^2} + C.
	\ee
\end{lem}
    
\begin{lem} For any smooth function $f$, we have
	\be\label{CJ}
	\int_\O\rho \f{D}{Dt} f dx = \f{d}{dt}\int_\O \rho f dx,
	\ee
where $\rho$ is the density satisfying $\eqref{eq1}_1$ and  $\frac{D}{Dt}f$ denotes the material derivative of a function f, i.e.
\be
\frac{D}{Dt}f=\p_t f+\u\cdot\nabla f.
\ee
\end{lem}

The following Zlotnik inequality will be used to get the time-independent upper bound of the
density.
\begin{lem}\label{zlemma} (\cite{zlot}) Let the function $y\in W^{1,1}(0, T)$ satisfys
\[
y'(t) = g(y) + h'(t) \text{ on } [0, T],\quad y(0) = y_0,
\]
with $g\in C(R)$ and $h\in W^{1,1}(0, T)$. If $g(\infty)=-\infty$ and
\be
h(t_2) - h(t_1) \leq N_0 + N_1(t_2 - t_1)
\ee
for all $0 \le t_1 < t_2 \le T$ with some $N_0 \ge 0$ and $N_1 \ge 0$, then
\[
y(t) \le \max(y_0, \zeta_0) + N_0 < \infty \text{ on } [0, T],
\]
where $\zeta$ is a constant such that $g(\zeta) \leq -N_1$ for $\zeta\ge\zeta_0$.
\end{lem}

\section{A priori estimates (I): upper bound of the density}

Let $T>0$ be a fixed time and $(\rho,\u,\B)$ with
$$\rho(x,t)=\rho(r,t),\ \u(x,t)=u(r,t)\dfrac{x}{r},\ \B(x,t)=B(r,t)\dfrac{x^{\perp}}{r}$$ be a strong solution to the problem \eqref{eq1}--\eqref{BD} on $\O\times(0, T)$. Then we will give some necessary a priori estimates for $(\rho,\u, \B)$. We begin with the following basic energy estimates.
\begin{lem}\label{EE} {\bf (Basic estimates)}
	There is a constant $C>0$, depending only on the initial data $(\rho_{0},  \u_{0}, \B_0$, such that
	\begin{equation}\label{EE2}
		\ba
		\sup_{t \in [0,\ T]} &(\| \sqrt{\rho }\u\|_{2}^{2} + \| \rho \|_{\gamma} ^{\gamma} + \|\B\|_2^2) \\
		&\quad + \int _{0}^{T}(\mu\| \nabla \u\|_{2}^{2} + \|(\mu+\lambda(\rho))^{\frac{1}{2}} \div \u\|_{2}^{2} +\|\g\B\|_2^2)  dt \leqs E_0.
		\ea
	\end{equation}
and for any $1\leq q<\infty$, it holds that
	\begin{equation}\label{BE2}
	\ba
	\sup_{t \in [0,\ T]} \|\B\|_q^q + \int_0^T \||\B|^{\f{q-2}{2}}|\g \B|\|_2^2 dt\leqs E_0.
	\ea
\end{equation}
and
	\begin{equation}\label{GBE2}
	\ba
	\sup_{t \in [0,\ T]} \|\g\B\|_2^2 + \int_0^T \|\lap \B\|_2^2 dt\leqs E_0.
	\ea
\end{equation}
where
\[
E_0 = \int_\O \rho_o|u_0|^2 + \rho_0^\gamma + |\B_0|^2 + |\g\B_0|^2dx.
\]
\end{lem}
\proof
\def\A{\mathbf{S}}
Denote that $\A=\u\cdot\g\B - \B\cdot\g\u + \B\div\u$, then
\[
%\int \lap\B\cdot \B_t = -\int \p_t\p_i B_j \p_i B_j dx = -\f 12\f{d}{dt}\int |\g\B|^2dx = \int \lap\B\cdot \lap\B dx -\int %\lap\B\cdot \A dx.
\int \lap\B\cdot \B_t = -\f 12\f{d}{dt}\int |\g\B|^2dx = \int \lap\B\cdot \lap\B dx -\int \lap\B\cdot \A dx.
\]
Therefore, by using lemma \ref{lem2}, one obtains that
\[\ba
\f 12\f{d}{dt}&\int |\g\B|^2dx + \int |\lap\B|^2 dx= -\int \lap\B\cdot \A dx\\
&\leqs \eps\|\lap\B\|_2^2 + \|\A\|^2\\
&\leqs \eps\|\lap\B\|_2^2 +\int |\u|^2|\g\B|^2 + |\B|^2|\g\u|^2 dx\\
&\leqs \eps\|\lap\B\|_2^2 +\|\g\u\|_2^2\int |\g\B|^2dx + \|\g\B\|_2^2\int|\g\u|^2 dx\\
&\leqs \eps\|\lap\B\|_2^2 +\|\g\u\|_2^2 \|\g\B\|_2^2.
\ea\]
\endproof

Define
~\be
\theta (\rho)= 2\mu \ln \rho+\frac{1}{\beta} (\rho^{\beta}  -1),
\ee
and the effective viscous flux
~\be\label{EVis}
F = (2\mu+\ld)\div \u - (P-\bar P) -\f 12 |\B|^2.
\ee
with $\bar P=\f{1}{|\O|}\int_\O Pdx$

By making use of the original mass equation $\eqref{eq2}_1$, we have
~\be\label{EDW}
\p_t\theta + u\partial_{r} \theta + F + P-\bar P +\f 12 {B^2} = 0.
\ee
Next, thanks to the definition of the effective viscous flux~\eqref{EVis}, we thus rewrite the momentum equation~$\eqref{eq2}_2$ as:
\be\label{XM}
\p_t(\rho u ) + \partial_{r} (\rho u^2) + \frac{ 1}{r} \rho u^2 = \partial_{r}F -\f{B^2}{r}.
\ee
Integrating over~$(R, r)$ gives
$$\frac{d}{dt} \int_{R}^{r} \rho u ds + \rho u^2 +\int_{R}^{r}\Big(\frac{\rho u^2}{s}+\f{B^2}{s}\Big) ds = F(r,t)- F(R,\ t).$$
Set
~\be\label{xieta}
\xi=\int_R^r\rho u dr,\quad \eta = \rho u^2 + \int_R^r\Big(\f{\rho u^2}{s}+\f{B^2}{s}\Big)ds.
\ee
One has
$$\xi_{t} + \eta - F + F(R,\ t)=0.$$
Adding the above equation to \eqref{EDW} gives
~\be\label{thxi}
(\theta +\xi )_{t} + u \cdot \partial_{r}(\theta+\xi) +\int_{R}^{r}\Big( \frac{\rho u^2}{s} +\f{B^2}{s} \Big)ds +P-\bar P + \f{B^2}{2}+ F(R,\ t)=0.
\ee

Indeed, \eqref{thxi} is exact  the corresponding form as (3.42) in \cite{2016Huang} in radially symmetric case which is the key commutator's estimates combined with the density equation and variable coefficients.

Note that $F(R,t) \neq 0$ due to the Dirichlet boundary condition \eqref{BD}. Then we will adopt the proof of  Lemma 3.2 first used in our previous paper \cite{2023Huang}. 

The key observation is that we could not get the estimate $\|F(R,t)\|_{L^1(0,T)}$ directly but instead to bound $\Big|\int_{0}^{t} F(R, s) ds \Big|$.

\begin{lemma} \label{FR}
It holds that
	\be\label{frrr}
	F(R,\ t) = \frac{1}{R^2} \Big[ \frac{d}{dt} \int_{0}^{R} \rho u r^2 dr + 2\int_{0}^{R} F(r, t) r dr - \int_{0}^{R} (\rho u^2 -B^2)r ds\Big],
	\ee
	and, for all $0\leq t\leq T$,
	\be
	\Big|\int_{0}^{t} F(R, s) ds \Big|\leqs 1 + \int_0^t\|\rho\|_{\beta}^{\beta}d\ts.
	\ee
\end{lemma}

\proof Multiplying~\eqref{XM} by~$r^2$ gives
$$\p_t(r^2 \rho u)+\partial_{r}(r^2 \rho u^2) - r (\rho u ^2 - B^2)= r^2\partial _{r}F,$$
and then integrating over~$(0, R)$ leads to
$$\int_{0}^{R}\p_t(r^2\rho u) dr + \int_{0}^{R} \partial_{r}(r^2\rho u^2) dr -\int_{0}^{R} r (\rho u^2 - B^2)dr=\int_{0}^{R} r^2\partial_{r}Fdr.$$
Observing that
$$\int_{0}^{R} \partial_{r}(r^2\rho u^2) dr =r^2\rho u^2 \big|_{0}^{R} = 0,$$
and 
$$\int_{0}^{R} r^2\partial_{r}Fdr = r^2 F \big |_{0}^{R}-\int_{0}^{R} F\cdot 2rdr =R^2F(R,\ t)- 2\int_{0}^{R} F\cdot rdr,$$
one has
$$F(R,\ t)= \frac{1}{R^2}\Big [ \frac{d}{dt}\int_{0}^{R} r^2 \rho u dr+ 2\int_{0}^{R} F\cdot rdr-\int_{0}^{R} r (\rho u^2-B^2) dr\Big].$$

Next, we come to estimate $F(R,t)$.\\
Integrating $F(R,t)$ over~$(0,\ t)$, we have
\begin{eqnarray*}
	\int_{0}^{t} F(R, \ts) d\ts &=&  \frac{1}{R^2} \Big[\int_{0}^{R} \rho u r^2 dr \big|_{\ts=t} - \int_{0}^{R} \rho u r^2 dr \big|_{\ts=0} + 2\int_{0}^{t}\int_{0}^{R} F(r, \ts) r dr d\ts\\
	& - &\int_{0}^{t}\int_{0}^{R} (\rho u^2-B^2) r dr d\ts \Big],
\end{eqnarray*}
and then
\begin{eqnarray*}\label{iFR}
	\begin{split}
		\Big|\int_{0}^{t} F(R,\ \ts) d\ts \Big| & \leq  \frac{1}{R^2} \Big[ R \cdot \int_{0}^{R} \rho |u| r dr + R \int_{0}^{R} \rho_{0}|u_{0}| r dr + 2\int_{0}^{t}\int_{0}^{R} F(r, \ts) r dr d\ts \\
		& \quad+  \int_{0}^{t}\int_{0}^{R} (\rho u^2 r +B^2) dr d\ts \Big]\\
		& \leqs 1 + \int_{0}^{T}\int_{0}^{R} |F(r, \ts)| r dr d\ts,
	\end{split}
\end{eqnarray*}
where we used the fact that~$|u| \leq 1+|u|^2$, and $\rho |u| r \leq \rho r + \rho |u|^2 r$.\\
It remains to estimate the term~$\int_{0}^{T}\int_{0}^{R} |F(r, \ts)| r dr d\ts $.
Notice that $$F = (2\mu + \lambda) \div u -(P-\bar P) - \f {1}2B^2,$$ so we have
\begin{eqnarray*}
	2\pi\int_{0}^{R} |F(r, \ts)| r dr &=& \int_{\Omega} |F(r, \ts)| dx  
	=  \int_{\Omega}|(2\mu + \lambda) \div u -(P-\bar P)- \f {1}2B^2| dx\\
	&  \leqs &\int_{\Omega}(2\mu + \lambda)|\div u| dx + 1\\
	& \leqs & \int_{\Omega} (2\mu + \lambda)|\div u|^2 dx + \|\rho\|_{\beta}^{\beta} +1.
\end{eqnarray*}

Hence we obtain that
\begin{eqnarray*}
	\Big|\int_{0}^{t} F(R,\ \ts) d\ts \Big| & \leqs & 1 + \int_{0}^{T}\int_{0}^{R} |F(r, \ts)| r dr d\ts \\
	& \leqs & 1+ \int_{0}^{T} \int_{\Omega} (2\mu + \lambda)|\div u|^{2} dx d\ts +\int_0^t\|\rho\|_{\beta}^{\beta}d\ts\\
	& \leqs & 1 + \int_0^t\|\rho\|_{\beta}^{\beta}d\ts,
\end{eqnarray*}
which completes the proof.
\endproof

%%c1

Next, we will give the High-order integrability of the density and velocity.\\
Let
~\be\label{Ge}
G(t)\triangleq\int_{0}^{t} F(R, s)ds + 2C_1 + C_2\int_0^T\|\rho\|_{\beta}^{\beta}d\ts,
\ee
with $C_1$ and $C_2$ large enough, one gets
$$1 \leq C_1 \leq G(t) \leq 3 C_1+ C_2\int_0^T\|\rho\|_{\beta}^{\beta}d\ts.$$
Then \eqref{thxi} could be rewritten as
~\be\label{thxiG}
(\theta + \xi + G(t))_{t} + u\cdot \nabla(\theta+\xi +G(t)) +\int_{R}^{r}\Big( \frac{\rho u^2}{s} + \f{B^2}{s} \Big)ds +P(r, t)+\f {B^2}2 =0.
\ee

Now, we are ready to prove the higher integrability of the density.
\begin{Proposition}\label{BJ}
	If~$\beta >1,\gamma>1$, then there exists a constant $C>0$ depending only on ~$\beta, \gamma, T$ and the initial data, such that
	\be
	\sup_{t\in[0,\ T]}\int_{\Omega} \rho^{2\beta\gamma+1 }dx \leq C.
	\ee
\end{Proposition}

\proof
Let~$f=(\theta + \xi + G(t))_{+}$,
multiplying the \eqref{thxiG} by ~$\rho f^{2\gamma -1}$,
and integrating over $\Omega$, we have
\begin{eqnarray*}
	\int_\O \rho f^{2\gamma -1}\cdot ((\theta &+& \xi + G(t))_{t} + u\cdot \nabla(\theta+\xi +G(t)))dx\\
	&=&  \frac{1}{2\gamma} \Big[\int_\O\rho \cdot\partial_{t}(f^{2\gamma} + \rho u\cdot \partial_{r}(f^{2\gamma})dx\Big]\\
	&=& \frac{1}{2\gamma} \Big[\int_\O\partial_{t}(\rho f^{2\gamma}) - \rho_{t}\cdot f^{2\gamma} + \div (\rho u\cdot f^{2\gamma}) - f^{2\gamma} \cdot \div(\rho u)dx\Big]\\
	&=&\frac{1}{2\gamma}\Big[\int_\O \partial_{t}(\rho f^{2\gamma}) + \div (\rho u\cdot f^{2\gamma})dx\Big].
\end{eqnarray*}
Notice that $$\int_{\O} \div (\rho u\cdot f^{2\gamma})dx = \int_{\partial \Omega} \rho u\cdot f^{2\gamma}\cdot n ds = 0.$$
We have
\[\ba
	\frac{1}{2\gamma} \frac{d}{dt}\int_{\O} &\rho f^{2\gamma} dx 
	\leqs \int_{\Omega} \Big[\int_{r}^{R} \Big(\frac{\rho u^2}{s} + \f{B^2}{s}\Big) ds +1\Big] (\rho f^{2\gamma-1})dx.\\
	& \leqs  \| \rho^{\frac{1}{2\gamma}} f \|_{L^{2\gamma}(\Omega)}^{2\gamma-1}\| \rho \|_{L^{2\beta\gamma+1}(\Omega)}^{\frac{1}{2 \gamma}} 
	\|\int_{r}^{R}  \Big(\frac{\rho u^2}{s} + \f{B^2}{s} \Big)ds\|_{L^{\frac{2\beta\gamma+1}{\beta}}(\Omega)}\\
	&\leqs  \| \rho^{\frac{1}{2\gamma}} f \|_{L^{2\gamma}(\Omega)}^{2\gamma-1} 
	\Big[\| \rho \|_{L^{2\beta\gamma+1}(\Omega)}^{\frac{1}{2 \gamma}} 
	\| \rho u^2 \frac{x}{r^2} + B^2 \frac{x}{r^2}\|_{L^{\frac{2(2\beta \gamma +1)}{2\beta\gamma+1+2\beta}}(\Omega)} +1\Big]\\
	&\leqs  \| \rho^{\frac{1}{2\gamma}} f \|_{L^{2\gamma}(\Omega)}^{2\gamma-1} 
	\Big[\| \rho \|_{L^{2\beta\gamma+1}(\Omega)}^{\frac{1}{2 \gamma}} 
	\Big(\| \rho u \|_{L^{\frac{2\beta\gamma+1}{\beta}}(\Omega)} \| \g u\|_{L^2(\Omega)} 
	+\|B\|_{L^{\frac{2\beta\gamma+1}{\beta}}(\Omega)}\|\g B\|_{L^2(\Omega)} \Big) +1\Big] ,
\ea\]
where we used Proposition~\ref{ME} and H\"{o}lder's inequality.\\
Since
\[\ba
\|\rho u \|_{L^{\frac{2\beta\gamma+1}{\beta}}(\Omega)} &\leqs \| \rho \|_{L^{2\beta\gamma+1}(\Omega)} \| u \|_{L^{\frac{2\beta\gamma+1}{\beta-1}}(\Omega)}\\
&\leqs\| \rho \|_{L^{2\beta\gamma+1}(\Omega)}\| \nabla u \|_{L^2(\Omega)},
\ea\]
we obtain that
~\be\label{fgam}
\frac{d}{dt}\int _{\Omega}\rho f^{2\gamma} dx  \leqs \Big(1+ \int_{\Omega}\rho f^{2\gamma } +\int _{\Omega}\rho^{2\beta\gamma+1 }dx\Big)
\Big( \| \g u \|_{L^{2}(\Omega)}^2
+\|\g B\|_2^2 +1\Big).
\ee

Next, we will estimate the term~$\int_{\Omega}\rho^{2\beta\gamma+1} dx$.\\
Notice that $f=(\theta(\rho) + \xi +G(t))_{+} \geq (\theta(\rho) + \xi +G(t))$ and
$$\rho^{\beta} \leqs \theta(\rho),\text{ when }\rho\geq 2.$$
For any~$\rho\geq 2 $, we have
$$\rho^{\beta}\leqs \theta(\rho) \leqs (f-\xi -G(t)) \leqs (f+|\xi|).$$
Hence, we can divide the term~$\int_{\Omega} \rho^{2\beta\gamma+1} dx$ into two parts:
\begin{eqnarray*}
	\int_{\Omega} \rho^{2\beta\gamma+1} dx &=&\int_{\Omega \cap {\rho (\leq 2)} } \rho^{2\beta\gamma+1} dx +\int_{\Omega \cap (\rho > 2) } \rho^{2\beta\gamma+1} dx\\
	& \leqs & 1+\int_{\Omega \cap (\rho >2)} \rho f^{2\gamma} dx +\int_{\Omega}\rho |\xi|^{2\gamma} dx\\
	& \leqs & 1+\int_{\Omega \cap (\rho >2)} \rho f^{2\gamma} dx  + \| \rho \|_{L^{\frac{2\beta\gamma+1}{2\beta\gamma+1-\gamma}}(\Omega)}  \| \xi \|_{L^{2(2\beta\gamma+1)}(\Omega)}^{2\gamma}\\
	& \leqs & 1+\int_{\Omega \cap (\rho >2)} \rho f^{2\gamma} dx  + \| \rho \|_{L^{\frac{2\beta\gamma+1}{2\beta\gamma+1-\gamma}}(\Omega)}  \| \nabla \xi \|_{L^{\frac{2\beta\gamma+1}{\beta\gamma+1}}(\Omega)}^{2\gamma}\\
	& \leqs & 1+\int_{\Omega \cap (\rho >2)} \rho f^{2\gamma} dx  + \| \rho \|_{L^{\frac{2\beta\gamma+1}{2\beta\gamma+1-\gamma}}(\Omega)} \|\rho^{\frac{1}{2}}\|_{L^{2(2\beta\gamma+1)}(\Omega)}^{2\gamma} \|\rho ^{\frac{1}{2}} u\|_{L^2(\Omega)}^{2\gamma},
\end{eqnarray*}
Since $$\frac{2\beta\gamma+1}{2\beta\gamma+1-\gamma} \leq 2\beta\gamma+1,$$
we have
\begin{eqnarray*}
	\| \rho \|_{L^{\frac{2\beta\gamma+1}{2\beta\gamma+1-\gamma}}(\Omega)} \|\rho^{\frac{1}{2}}\|_{L^{2(2\beta\gamma+1)}(\Omega)}^{2\gamma}
	&\leqs&\| \rho \|_{2\beta\gamma+1} \| \rho \|_{2\beta\gamma+1}^{\gamma}\leqs\| \rho \|_{2\beta\gamma+1}^{\gamma+1}.
\end{eqnarray*}
For~$\| \rho \|_{2\beta\gamma+1}^{\gamma+1}$,
taking~$p=\frac{2\beta\gamma+1}{\gamma+1}, \quad q=\frac{2\beta\gamma+1}{2\beta\gamma-\gamma}$ satisfying $\frac{1}{p} + \frac{1}{q}=1$, and using Young's inequality yield
\[\| \rho \|_{2\beta\gamma+1}^{\gamma+1} \leqs \eps (\| \rho \|_{2\beta\gamma+1}^{\gamma+1})^{\frac{2\beta\gamma+1}{\gamma+1}}+1.
\]
Take~$\eps$ small enough, we obtain that
\begin{eqnarray*}
	\int_{\Omega} \rho^{2\beta\gamma+1} dx &\leqs &1+ \int_{\Omega}\rho f^{2\gamma} dx +\| \rho \|_{2\beta\gamma+1}^{\gamma+1}\\
	&\leqs & 1+ \int_{\Omega}\rho f^{2\gamma} dx +\frac{1}{2}\| \rho \|_{2\beta\gamma+1}^{2\beta\gamma+1}.
\end{eqnarray*}
Therefore,
$$\int_{\Omega} \rho^{2\beta\gamma+1} dx \leqs 1+ \int_{\Omega} \rho f^{2\gamma} dx.$$
Submitting the above into~\eqref{fgam} gives
~\be
\frac{d}{dt}\int _{\Omega}\rho f^{2\gamma} dx  \leqs \Big(1+ \int_{\Omega}\rho f^{2\gamma }dx \Big) 
\Big(\| \g u \|_{L^{2}(\Omega)}^2 +\|\g B\|_2^2+1\Big).
\ee
Finally, ~Gronwall's inequality implies that
$$\sup_{t\in[0,\ T]}\int_{\Omega} \rho^{2\beta\gamma+1 }dx \leqs \sup_{t\in[0,\ T]}\int _{\Omega}\rho f^{2\gamma} dx +1 \leqs 1.$$
This completes the proof.
\endproof

\begin{Proposition}\label{QS} If
	~\be
	0<\v <\min\Big\{2,\ \f{4\gamma(\beta-1)+2}{2\gamma-1}\Big\},
	\ee
	then there exists a constant $C_\v>0$ depending on the time T such that
	~\be\label{eudelta}
	\sup_{t\in [0,\ T]}\int_{\Omega} \rho |u|^{\v+2} dx\leqs C_\v(T).
	\ee
\end{Proposition}

\proof
Multiplying~$\eqref{eq1}_{2}$ by~$(2+\v)u\cdot |u|^{\v}$, and integrating (by parts) over $\Omega$ gives  
~\be\label{Udelta}
\ba
\frac{d}{dt} \int_{\Omega}\rho|u|^{\v+2} dx &+ \mu(\v+2) \int_{\Omega}\Big[\v|u|^{\v}\cdot |\nabla |u||^2 + |u|^{\v} |\nabla u|^{2}\Big]dx\\
&+ (2+\v) \int_{ \Omega }(\mu+ \lambda(\rho)) |u|^{\v} |\div u|^{2}dx\\
&+ (2+\v) \int_{ \Omega }\v|u|^{\v-1}  (u \cdot \nabla) |u| \cdot (\mu+ \lambda(\rho))\div u dx\\
= &\ (2+\v)\int_{\Omega} (P+\f 12 B^2) \cdot \div (|u|^{\v}\cdot u) dx +(2+\v)\int_\O |u|^\v B\cdot\g B\cdot udx\\
\leqs& \ \int_{\Omega} (\rho^{\gamma} +1)|u|^{\v} |\nabla u| dx.
\ea\ee
Here we used the fact that
\begin{eqnarray*}
	\div (|u|^{\v} \cdot u) &=& u \cdot \nabla |u|^{\v} + |u|^{\v} \cdot \div u\\
	& \leqs & |u|^{\v}  |\nabla u|.
\end{eqnarray*}

Notice that
\[\ba
(2+\v)\v \int_{ \Omega } &(\mu+ \lambda(\rho))\div u |u|^{\v-1}  (u \cdot \nabla) |u|  dx\\
&=(2+\v)\v \int_{ \Omega } (\mu+ \lambda(\rho))|u|^{\v}\big(|\p_r u|^2 + \f u r\p_r u\big) dx\\
&\geq (2+\v)\v \int_{ \Omega } (\mu+ \lambda(\rho))|u|^{\v}\big(|\p_r u|^2 -\f 12|\div u|^2\big) dx,
\ea\]
and
\[\ba
(2+\v)\int_\O &|u|^\v B\cdot\g B\cdot udx
= -(2+\v)\int_\O (B\otimes B) :\g( |u|^\v u)dx \\
&\leqs \int_\O |B|^2 |u|^\v|\g u|dx,
\ea\]
which combined with~\eqref{Udelta} gives
\[\ba
\frac{d}{dt} \int_{\Omega}\rho|u|^{\v+2} dx &+ \mu(\v+2) \int_{\Omega}[\v|u|^{\v}\cdot |\nabla |u||^2 + |u|^{\v} |\nabla u|^{2}]dx\\
&+ (2+\v) \int_{ \Omega }(\mu+ \lambda(\rho)) |u|^{\v} |\div u|^{2}dx\\
&+ (2+\v)\v \int_{ \Omega } (\mu+ \lambda(\rho))|u|^{\v}\big(|\p_r u|^2 -\f 12|\div u|^2\big) dx\\
\leqs & \int_{\Omega} (\rho^{\gamma}+|B|^2) |u|^{\v} |\nabla u| dx.
\ea\]
Since~$1-\frac{\v}{2} > 0 $, after applying ~H\"{o}lder's inequality and ~Young's inequality, we obtain that
\[\ba
\frac{d}{dt}&\int_{\Omega}\rho |u|^{2+\v} dx + \mu(\v+2)\int_{\Omega}|u|^{\v} |\nabla u|^{2} dx
\leqs \int_{\Omega}(\rho^{\gamma}+|B|^2)|u|^{\v} |\nabla u| dx\\
& \leqs \eps\int_{\Omega} |u|^{\v}|\nabla u|^{2} dx + \int_{\Omega} \rho |u|^{2+\v} dx +
\int_{\Omega} \rho^{\gamma(2+\v)-\f \v 2} dx + \int_{\Omega} |B|^4|u|^\v dx\\
& \leqs \eps\int_{\Omega} |u|^{\v}|\nabla u|^{2} dx + \int_{\Omega} \rho |u|^{2+\v} dx +
\int_{\Omega} \rho^{\gamma(2+\v)-\f \v 2} dx + \int_{\Omega} |B|^4(1+|u|) dx\\
&\leqs \eps\int_{\Omega} |u|^{\v}|\nabla u|^{2} dx + \int_{\Omega} \rho |u|^{2+\v} dx +
\int_{\Omega} \rho^{\gamma(2+\v)-\f \v 2} dx + \|B\|_4^4 + \|B\|_3^3\|\g B\|_2\g u\|_2.
\ea\]
Hence, by taking $\eps$ small enough, \eqref{eudelta} follows from Proposition~\ref{BJ}, the fact  $\gamma(2+\v)-\f \v 2 \leq 2\beta\gamma+1$, and~Gronwall's inequality immediately.
\endproof

Next, we will give the upper bound of the density.
\begin{Proposition}\label{SZ} Under the conditions of Theorem \ref{them1}, there is a constant $C>0$ depending on the time T, such that
	~\be
	\| \rho \|_{L^{\infty} L^{\infty}}  \leq C.
	\ee
\end{Proposition}	 
%%c1

\proof Let ~$\t G(t)\triangleq\int_{0}^{t} F(R, \ts)d\ts$. Recalling the equation of $\theta$,
$$(\theta + \xi + \t G(t))_{t} + u\cdot \nabla(\theta+\xi +\t G(t)) +\int_{R}^{r} \Big( \frac{\rho u^2}{s}+\f{B^2}{s} \Big)ds +(P-\bar P)+\f{B^2}{2} =0.$$
we have
$$
\ba
(\theta + \xi + \t G(t))_{t} &+ u\cdot \nabla(\theta+\xi +\t G(t))  +P-\bar P+\f{B^2}2 = -\int_{R}^{r} \Big(\frac{\rho u^2}{s}+\f{B^2}{s} \Big) ds\\
&\leq \|\rho\|_{L^\infty L^\infty}\|\g u\|_2^2 + \|\g B\|_2^2 +1.
\ea$$
where we used the following fact
\[
2\pi\int_0^R \Big(|\p_r u|^2r + \f{u^2}{r}\Big)dr
=\|\div \u\|_2^2\leqs \|\g u\|_2^2.
\]

Integrating above inequality on ~$(0,\ t)$, we get
\begin{eqnarray*}
	(\theta + \xi + \t G(t))(r, t) & \leq &\|(\theta + \xi )(r, 0)\|_{\infty}+\|\rho\|_{L^\infty L^\infty}\int_{0}^{t} \|\g u\|_2^2 d\ts +\int_0^t\|\g B\|_2^2ds+1\\
	& \leqs & 1 + \| \rho \|_{L^{\infty} L^{\infty}}.
\end{eqnarray*}
By proposition ~\ref{QS} and ~H\"{o}lder's inequality,
\begin{eqnarray*}
	|\xi| \leq \int_{0}^{R} |\rho u| dr & = & \Big |\int_{0}^{R} \rho^{\frac{1}{2+\v}} \cdot \rho^{\frac{1+\v}{2+\v }} \cdot |u| \cdot r^{\frac{1}{2+\v}} \cdot r^{-\frac{1}{2+\v}} dr \Big|\\
	& \leq &\Big [\int_{0}^{R}(\rho^{\frac{1}{2+\v}} \cdot |u| \cdot r^{\frac{1}{2+\v}} )^{2+\v}dr  \Big]^{\frac{1}{2+\v}} \cdot
	\Big[ \int_{0}^{R} (\rho^{\frac{1+\v}{2+\v}}r^{-\frac{1}{2+\v}} )^{\frac{\v+2}{\v+1}} dr \Big]^{\frac{\v+1}{\v+2}}\\
	&\leqs & C(T)\| \rho \|_{L^{\infty}}^{\frac{1+\v}{2+\v}} (\int_{0}^{R} r^{-\frac{1}{1+\v}} dr)^{\frac{\v+1}{\v+2}}\\
	& \leqs &C(T) \| \rho \|_{L^{\infty}}^{\frac{1+\v}{2+\v}}.
\end{eqnarray*}
Applying lemma ~\ref{FR}, it holds that
\[
|\t G(t)|\leqs 1+ \int_0^t\|\rho\|_\beta^{\beta} d\ts\leqs 1+ \int_0^t\|\rho\|_\infty^{\beta} d\ts
\]
Therefore, we get
\begin{eqnarray*}
	\theta &\leq& (\theta+\xi+\t G) +|\xi| + |\t G|\\
	& \leqs &  1+ \| \rho \|_{L^{\infty} L^{\infty}} +\| \rho \|_{L^{\infty}L^{\infty}}^{\frac{1+\v}{2+\v}}+ \int_0^t\|\rho\|_\infty^{\beta} d\ts\\
	& \leqs &  1+ \| \rho \|_{L^{\infty} L^{\infty}}+ \int_0^t\|\rho\|_\infty^{\beta} d\ts).
\end{eqnarray*}

Notice that ~$\theta(\rho)= 2\mu \ln \rho + \frac{1}{\beta}(\rho^{\beta}-1)$, we have
$$\ba
\rho^{\beta} &\leq \max (2^{\beta}, \theta(\rho))\\
&\leqs 1+ \| \rho \|_{L^{\infty} L^{\infty}}+ \int_0^t\|\rho\|_\infty^{\beta} d\ts.\\
\ea$$
Therefore, 
~\be\label{SQJ}
\|\rho\|_\infty^{\beta}\leqs 1+ \| \rho \|_{L^{\infty} L^{\infty}}+ \int_0^t\|\rho\|_\infty^{\beta} d\ts.
\ee
Applying ~Gronwall's inequality, it follows that
\[
\int_0^t\|\rho\|_\infty^{\beta} d\ts\leqs 1+ \| \rho \|_{L^{\infty} L^{\infty}}.
\]
This, together with~\eqref{SQJ}, gives
\[
\|\rho\|_\infty^{\beta}\leqs 1+ \| \rho \|_{L^{\infty} L^{\infty}}.
\]
Taking superium according to ~$t\in [0,T]$, one gets
\[
\|\rho\|_{L^{\infty}L^{\infty}}^{\beta}\leqs 1+ \| \rho \|_{L^{\infty} L^{\infty}}.
\]
Applying ~Young inequality, we have
\[
\|\rho\|_{L^{\infty}L^{\infty}}^{\beta}\leqs \f 12 \|\rho\|_{L^\infty L^\infty}^{\beta} + 1.
\]
This gives
$$\| \rho \|_{L^{\infty}L^{\infty}}^{\beta} \leq C(T),$$
and finishes the proof of Proposition \ref{SZ}.
		  
%%c2
\endproof

\section{Upper bound of the density independent of time}
In this section, we denote $R_T=1+\sup\limits_{0<t<T}\r$. 

Here we will borrow the steps adapted in the proof of \cite{2016Huang} and carefully deal with extra terms caused by the commutator's estimate \eqref{thxi}.

The first Lemma is the corresponding version of Lemma 3.1 in  \cite{2016Huang} which can be much simplified for the radially symmetric case. One can also refer to \cite{2023Huang} for the compressible Navier-Stokes equations.

\begin{lemma} {\rm (Estimate of $F$)} For any $\alpha\in (0,1)$, there is a constant $C(\alpha)$ depending only on $\alpha, \mu,\beta,\gamma, \|\rho_0\|_{L^\infty}$, and $\|u_0\|_{H^1}$ such that
	\be\label{A2}
	\sup\limits_{0\leq t\leq T}\log\Big(e + \int\f{F^2}{2\mu+\ld} dx + \int(2\mu+\ld)|\div u|^2dx\Big) \leqs C(\a)R_T^{\a\beta+1} + R_T^{\gamma-2\beta+ 1}.
	\ee
\end{lemma}
\proof
The momentum equations can be rewritten as
\[
\r\dot u = \g \tF + B\cdot\g B.
\]
Consequently, multiplying the above equations by $\dot{u}$ and integrating the resulting equality over the ball $\Omega$ yields that
\[
\int \r |\dot u|^2 dx = \int \g \tF\cdot \dot u dx +\int B\cdot\g B\cdot\dot udx = -\int \tF\div \dot u dx+\int B\cdot\g B\cdot\dot udx.
\]
Note that
\[
\div \dot u = \f{D}{Dt} \div u + |\div u|^2 - 2\p_r u\f{u}{r}.
\]
Hence,
\[
\int \tF\div \dot u dx = \int \tF \f{D}{Dt} \div udx + \int \tF|\div u|^2dx - 2\int \tF\p_r u\f{u}{r}dx
\]
The first term can be estimated as follows
\[\ba
\int &\tF \f{D}{Dt} \div u dx  =\int \tF\f{D}{Dt} (\f{\tF+P-\bar P}{2\mu+\ld})dx+ \f12\int \tF D_t\f{B^2}{2\mu+\ld}\\
&=\f 12\f{d}{dt}\int\f{\tF^2}{2\mu+\ld} dx -\f 12\int \f{\tF^2}{2\mu+\ld}\div u dx\\
&\qquad + \f 12\int \tF^2 \f{\r\ld'}{(2\mu+\ld)^2}\div udx  - \int \tF  \r(\f{P}{2\mu+\ld})'\div udx\\
&\qquad + (\gamma-1)\int \f{F}{2\mu+\ld}dx\int(P-\bar P)\div u dx\\
&\qquad + \f12\int \tF\f{B^2\rho\ld'}{(2\mu+\ld)^2}\div u dx +\int \f{\tF BD_t B}{2\mu+\ld}dx
\ea\]
\iffalse
Therefore,
\[\ba
\int \r |\dot u|^2 dx & = -\int \tF\div \dot u dx+\int B\cdot\g B\cdot\dot udx\\
&=-\int \tF \f{D}{Dt} \div udx - \int \tF|\div u|^2dx + 2\int \tF\p_r u\f{u}{r}dx+\int B\cdot\g B\cdot\dot udx\\
&=-\f 12\f{d}{dt}\int\f{\tF^2}{2\mu+\ld} dx +\f 12\int \f{\tF^2}{2\mu+\ld}\div u dx- \f 12\int \tF^2 \f{\r\ld'}{(2\mu+\ld)^2}\div udx \\
&\qquad  + \int \tF  \r(\f{P}{2\mu+\ld})'\div udx + (\gamma-1)\int \f{F}{2\mu+\ld}dx\int(P-\bar P)\div u dx\\
&\qquad - \f 12\int \tF\f{B^2\rho\ld'}{(2\mu+\ld)^2}\div u dx -\int \f{\tF BD_t B}{2\mu+\ld}dx\\
&\qquad - \int \tF|\div u|^2dx + 2\int \tF\p_r u\f{u}{r}dx +\int B\cdot\g B\cdot\dot udx\\
\ea\]
\fi
Combing all the above inequalities together to get
\[\ba
\f 12\f{d}{dt}&\int\f{\tF^2}{2\mu+\ld} dx + \int \r |\dot u|^2 dx = 
\f 12\int \f{\tF^2}{2\mu+\ld}\div u dx- \f 12\int \tF^2 \f{\r\ld'}{(2\mu+\ld)^2}\div udx \\
&\qquad  + \int \tF  \r(\f{P}{2\mu+\ld})'\div udx + (\gamma-1)\int \f{F}{2\mu+\ld}dx\int(P-\bar P)\div u dx\\
&\qquad - \f 12\int \tF\f{B^2\rho\ld'}{(2\mu+\ld)^2}\div u dx -\int \f{\tF BD_t B}{2\mu+\ld}dx\\
&\qquad - \int \tF|\div u|^2dx + 2\int \tF\p_r u\f{u}{r}dx +\int B\cdot\g B\cdot\dot udx\\
&\leqs \int  \f{\tF^2}{(2\mu+\ld)}|\div u|dx  + \int (\f{|\tF|P}{2\mu+\ld})|\div u|dx
+ \|\g \tF\|_2\|\div u\|_2^2\\
&\qquad + (\gamma-1)\int \f{F}{2\mu+\ld}dx\int(P-\bar P)\div u dx +\int \f{|F|}{2\mu+\ld}|\div u|dx\\
&\qquad  + \int \f{|\tF|B^2}{2\mu+\ld}|\div u| dx +\int \f{\tF BD_t B}{2\mu+\ld}dx +\int B\cdot\g B\cdot\dot udx\\
&\leqs \int  \f{\tF^2}{(2\mu+\ld)}|\div u|dx  + \int (\f{|\tF|P}{2\mu+\ld})|\div u|dx  + \|\g \tF\|_2\|\div u\|_2^2\\
&\qquad + (\gamma-1)\int \f{F}{2\mu+\ld}dx\int(P-\bar P)\div u dx +\int \f{|F|}{2\mu+\ld}|\div u|dx\\
&\qquad + \int \f{|\tF|B^2}{2\mu+\ld}|\div u| dx +\int \f{\tF BD_t B}{2\mu+\ld}dx +\int B\cdot\g B\cdot\dot udx\\
&= I_1 + I_2 + I_3 + I_4 + I_5 + I_6 + I_7 + I_8.
\ea\]
Here $I_6-I_8$ are the new terms from the magnetic field which is different from the compressible Navier-Stokes equations.

We shall first bound $F$ which will play a crucial role in deriving the estimates for each $I_i$.

\[\ba
\|\tF\|_{H^1} &\leqs \Big|\int \tF dx\Big| + \|\g \tF\|_2\\
&\leqs \|\r\|_\beta^{\f{\beta}{2}}\|\sqrt{2\mu+\ld}\div u\|_2 + \|B\|_2^2+ R_T^{\f12}\|\sqrt\r\dot u\|_2 + \||B||\g B|\|_2.
\ea\]

Then we will estimate each $I_i$ as follows.

For any $0<\a<1$,
\[\ba
I_1 &=\int  \f{\tF^2}{(2\mu+\ld)}|\div u|dx \leqs \|\f{\tF^2}{2\mu+\ld}\|_2\|\div u\|_2\\
&\leqs \|\f{\tF}{\sqrt{2\mu+\ld}}\|_2^{1-\a} \|\tF\|^{1+\a}_{\f{2(1+\a)}{\a}}\|\div u\|_2\\
&\leqs \|\f{\tF}{\sqrt{2\mu+\ld}}\|_2^{1-\a} \|\tF\|_2^\a \|\tF\|_{H^1}\|\div u\|_2\\
&\leqs R_T^{\f {\a\beta}2}\|\f{\tF}{\sqrt{2\mu+\ld}}\|_2  \|\div u\|_2\|\tF\|_{H^1}\\
&\leqs \eps\|\sqrt{\r} \dot u\|_2^2 +\|\sqrt{2\mu+\ld}\div u\|_2^2 +\|\g B\|_2^2 + \||B||\g B|\|_2^2\\
&\qquad + C(\a)R_T^{\a\beta+1}\|\f{\tF}{\sqrt{2\mu+\ld}}\|_2^2\|\div u\|_2^2.
\ea\]
Then we will deal with $I_2$.
\[\ba
I_2 &=\int \f{|\tF|P}{2\mu+\ld}|\div u|dx \leq \|\tF\|_{2+4\gamma/\beta}\|\f{P}{(2\mu+\ld)^{\f 32}}\|_{2+\beta/\gamma}\|\sqrt{2\mu+\ld}\div u\|_2\\
&\leqs R_T^{\gamma/2-\beta}\|\tF\|_{H^1}\|\sqrt{2\mu+\ld}\div u\|_2\\
&\leqs \eps\|\sqrt{\r} \dot u\|_2^2  +\|\g B\|_2^2 + \||B||\g B|\|_2^2 
+ R_T^{\gamma-2\beta+ 1}\|\sqrt{2\mu+\ld}\div u\|_2^2.
\ea\]\\
Then we are about to estimate $I_3$.
\[\ba
I_3 &= \|\g \tF\|_2\|\div u\|_2^2 =\|\g F\|_2\|\f{\tF+(P-\bar P)+\f 12B^2}{2\mu+\ld}\|_2\|\div u\|_2\\
&\leqs \eps \|\sqrt{\r} \dot u\|_2^2 +\||B||\g B|\|_2^2+ R_T\|\f{\tF+P-\bar P+\f 12B^2}{2\mu+\ld}\|_2^2\|\div u\|_2^2\\
&\leqs \eps \|\sqrt{\r} \dot u\|_2^2 +\||B||\g B|\|_2+ R_T\Big(\|\f{\tF}{2\mu+\ld}\|_2^2+ \|\f{P-\bar P}{2\mu+\ld}\|_2^2 +1\Big)\|\div u\|_2^2\\
&\leqs \eps \|\sqrt{\r} \dot u\|_2^2 +\||B||\g B|\|_2^2+ R_T\|\f{\tF}{\sqrt{2\mu+\ld}}\|_2^2\|\div u\|_2^2+ R_T(\|\f{P-\bar P}{2\mu+\ld}\|_2^2+1)\|\div u\|_2^2\\
&\leqs \eps \|\sqrt{\r} \dot u\|_2^2 +\||B||\g B|\|_2^2+ R_T\|\f{\tF}{\sqrt{2\mu+\ld}}\|_2^2\|\div u\|_2^2+ R_T^{\gamma-2\beta+1}\|\div u\|_2^2.
\ea\]

$I_4$ can be bounded as follows
\[\ba
I_4 &\leqs \Big|\int \f{F}{2\mu+\ld}dx\int(P-\bar P)\div u dx\Big|\\
&\leqs \int|P-\bar P||\div u| dx\\
&\leqs \int|F-(2\mu+\ld)\div u -\f 12 B^2||\div u| dx\\
&\leqs \|F\|_{H^1}\|\div u\|_2 + \|\sqrt{2\mu+\ld}\div u\|_2^2 + \|\g B\|\|\g u\|.\\
&\leqs \eps \|\sqrt{\r}\dot u\|_2^2 + R_T\Big(\|\sqrt{2\mu+\ld}\div u\|_2^2 +\|\g B\|_2^2 
+\|\g u\|_2^2+\||B||\g B|\|_2^2\Big).
\ea\]

The estimate of $I_5$ is same as $I_4$.

To estimate $I_6$, by lemma \ref{lem2} and H\"older inequality, we have
\[\ba
I_6&\leqs \int |\tF| B^2|\g u| dx \leqs \|\tF\|_{H^1}\|\g u\|_2 \\
&\leqs (\|\r\|_\beta^{\f{\beta}{2}}\|\sqrt{2\mu+\ld}\div u\|_2 + \|B\|_2^2+ R_T^{\f12}\|\sqrt\r\dot u\|_2 + \||B||\g B|\|_2)\|\g u\|_2\\
&\leqs  \eps\|\sqrt\r \dot u\|_2^2 +\|\sqrt{2\mu+\ld}\div u\|_2^2 +\|\g B\|_2^2 + \||B||\g B|\|_2^2 + R_T\|\g u\|_2^2.
\ea\]

Now, using the equation of $\B$ and the estimate $I_7$, it follows that
\be\ba\label{I7}
I_7&=\int \f{\tF BD_t B}{2\mu+\ld}dx = \int \f{\tF B}{2\mu+\ld}(B\cdot\g u -  B\div u +\lap B)dx\\
&\leqs \int |\tF| B^2|\g u| dx + \int \f{\tF B}{2\mu+\ld}\lap Bdx.
\ea\ee
The first term of right hand side of \eqref{I7} can be estimated as $I_6$.
 The last term of above inequality can be estimated as following:
\[\ba
\int &\f{\tF B}{2\mu+\ld}\lap Bdx\leqs \|\f{\tF}{\sqrt{2\mu+\ld}}\|_2^{\f{1-\a}{2}} \Big(\int|\tF|^2|B|^{\f{4}{1+\a}}dx\Big)^{\f{1+\a}{4}} \|\lap B\|_2\\
&\leqs \|\f{\tF}{\sqrt{2\mu+\ld}}\|_2^{\f{1-\a}{2}} \|\tF\|^{\f{1+\a}{2}}_4\|B\|_{\f{8}{1+\a}}\|\lap B\|_2\\
&\leqs C(\a)R_T^{\f{\a\beta}{4}}\|\f{\tF}{\sqrt{2\mu+\ld}}\|_2^{\f{1}{2}} \|\tF\|_{H^1}^{\f{1}{2}} \|\lap B\|_2 \|\g B\|_2^{\f 12}\\
&\leqs \|\lap B\|_2^2 + C(\a)R_T^{\f{\a\beta}{2}}\|\f{\tF}{\sqrt{2\mu+\ld}}\|_2 \|\g B\|_2 \|\tF\|_{H^1}\\
&\leqs \eps\|\sqrt{\r} \dot u\|_2^2 +\|\sqrt{2\mu+\ld}\div u\|_2^2 +\|\g B\|_2^2 + \||B||\g B|\|_2^2+  \|\lap B\|_2^2\\
&\qquad + C(\a)R_T^{\a\beta + 1}\|\f{\tF}{\sqrt{2\mu+\ld}}\|_2^2\|\g B\|_2^2.
\ea\]

Now, we deal with $I_8$.
\[\ba
I_8&=\int B\cdot\g B\cdot \dot u dx =\int B\cdot\g B\cdot \p_t u dx + \int B\cdot\g B\cdot (u\cdot\g)u dx\\
&=\int \div(B\otimes B)\cdot \p_t u dx + \int B\cdot\g B\cdot (u\cdot\g)u dx\\
&=-\int (B\otimes B):\p_t \g u dx + \int B\cdot\g B\cdot (u\cdot\g)u dx\\
&=-\f{d}{dt}\int (B\otimes B): \g u dx + \int (B_t\otimes B + B\otimes B_t)\cdot \g u dx\\
&\qquad + \int B\cdot\g B\cdot (u\cdot\g)u dx\\
&=-\f{d}{dt}I_{81} + I_{82} + I_{83}.
\ea\]
Then, $I_{81}, I_{82}, I_{83}$ can be estimated as following.
\[\ba
I_{81} &= \int (B\otimes B): \g u dx \leqs \|\g u\|_2\|B\|_4^2 \\
&\leqs \eps\|\div u\|_2^2 + 1.
\ea\]

\[\ba
I_{82} &=\int (B_t\otimes B + B\otimes B_t): \g u dx\\
&\leqs \int |u\cdot\g B - B\cdot\g u + B\div u - \lap B||B||\g u|dx\\
&\leqs \|\g u\|_2^2\|\g B\|_2^2  + \|B\|_\infty \|\lap B\|_2\|\g u\|_2\\
&\leqs \|\g u\|_2^2\|\g B\|_2^2  + \|\g B\|_2 \|\lap B\|_2\|\g u\|_2\\
&\leqs \|\lap B\|_2^2 + \|\g u\|_2^2. 
\ea\]

\[\ba
I_{83} &=\int B\cdot\g B\cdot (u\cdot\g)u dx\\
&\leqs \|\g u\|_2^2\|\g B\|_2^2
\leqs \|\g u\|_2^2.
\ea\]

Substitute all of estimates of $I_1-I_8$, we get
\[\ba
\f{d}{dt} &\int\f{\tF^2}{2\mu+\ld} + (B\otimes B): \g u dx + \int \r |\dot u|^2 dx
\leqs   \eps\|\sqrt{\r} \dot u\|_2^2+  \|\lap B\|_2^2 \\
& + R_T\Big(\|\g B\|_2^2 +\|\g u\|_2^2+\||B||\g B|\|_2^2\Big)\\
 &  + C(\a)R_T^{\a\beta+1}\|\f{\tF}{\sqrt{2\mu+\ld}}\|_2^2\Big(\|\div u\|_2^2+\|\g B\|_2^2\Big)\\
 & + R_T^{\gamma-2\beta+ 1}\|\sqrt{2\mu+\ld}\div u\|_2^2.
\ea\]

Taking $\eps$ small enough and by Gronwall's inequality, we finally arrive at
\[
\sup\limits_{0\leq t\leq T}\log(e + \int\f{\tF^2}{2\mu+\ld} dx) \leqs C(\a)R_T^{\a\beta+1} + R_T^{\gamma-2\beta+ 1}.
\]
and moreover,
\[\ba
\int (2\mu+\ld)|\div u|^2 dx &=\int (2\mu+\ld)|\f{\tF+P+\f 12B^2}{2\mu+\ld}|^2 dx\\
&\leqs \int \f{\tF^2}{2\mu+\ld}dx + \int \f{P^2}{2\mu+\ld}dx + 1\\
&\leqs \int \f{\tF^2}{2\mu+\ld}dx + R_T^{\gamma-\beta}\\
\ea\]
This gives the desired estimates of \eqref{A2} and thus finished the proof Lemma 6.1.
\endproof

Another crucial estimate lies on the $\xi$.
\begin{lemma} {\rm (Estimate of $\xi$)} For any $\alpha\in (0,1)$, there is a constant $C(\alpha)$ depending only on $\alpha, \mu,\beta,\gamma, \|\rho_0\|_{L^\infty}$, and $\|u_0\|_{H^1}$ such that
\[
\|\xi\|_\infty\le C(\a)R_T^{1+\a\beta} + R_T^{1+\gamma/2-\beta}.
\]
\end{lemma}

\proof
Indeed, by Lemma 2.5 , we obtain that for any $\alpha\in (0,1)$
	\[\ba
	\|\xi\|_\infty &\leqs \|\g\xi\|_2\log(e+\|\g\xi\|_3) + \|\xi\|_2 +1\\
	&\leqs \|\r u\|_2\log(e+\|\r u\|_3) + \|\r u\|_{\f{2\gamma}{\gamma+1}} +1\\
	&\leqs R_T^{1/2}\log(e+R_T\|\g u\|_2)  +1\\
	&\leqs R_T^{1/2}\log(e+\|\g u\|_2)  + R_T\\
	&\leqs R_T^{1/2}\Big(C(\a)R_T^{(1+\a\beta)/2} + R_T^{(\gamma-2\beta+1)/2}\Big)  +1\\
	&\leqs C(\a)R_T^{1+\a\beta} + R_T^{1+\gamma/2-\beta}\\
	\ea\]

\endproof

Next, we are ready to give the upper bound of the density independent of the time T.
\begin{Proposition}\label{SZ1} Under the conditions of Theorem \ref{them1}, and 
\be
\max\{1,\frac{\gamma+2}{4}\}<\beta\le\gamma.
\ee
then there is a constant $C>0$ independent of the time T, such that
	~\be
	\| \rho \|_{L^{\infty} L^{\infty}}  \leq C.
	\ee
\end{Proposition}
\proof
Recall \eqref{thxi} that
\be\label{thxi-1}
(\theta +\xi )_{t} + u \cdot \partial_{r}(\theta+\xi) +\int_{R}^{r}\Big( \frac{\rho u^2}{s} +\f{B^2}{s}\Big) ds +P-\bar P +\f{B^2}{2} + F(R,\ t)=0.
\ee

The main idea is to employ Zlotnik's inequality to get the time-independent estimate for the density. We just need to rewrite the equation in suitable form and verify the conditions required by the Zloinik type inequality. This will result in the restriction of $\beta$.

Set 
\be
y=\theta+\xi,\ g(y)=-P(\theta^{-1}(y-\xi)),\ g(+\infty)=-\infty,
\ee
and
\be
h=\int_0^t\int_R^r \Big(\f{\r u^2}{s}+\f{B^2}{s}\Big)dsd\tau + \int_0^t\Big( F(R,\tau) - \bar P + \f{B^2}{2}\Big)d\tau.
\ee
Then \eqref{thxi-1} can be written as
\be
\frac{D}{Dt}y=g(y)+h_t.
\ee
It suffices to verifies the conditions required in lemma \ref{zlemma}, especially the Lipschitz coefficients of $h$.

For any time $0<t_1<t_2$,
\[
h(t_2)-h(t_1)=\int_{t_1}^{t_2}\int_R^r \Big(\f{\r u^2}{s}+\f{B^2}{s}\Big)dsd\tau 
+ \int_{t_1}^{t_2} \Big( F(R,\tau) - \bar P(\tau)+\f{B^2}{2}\Big)d\tau
\]

\[\ba
|h(t_2)-h(t_1)|&=\int\int_R^r \Big(\f{\r u^2}{s}+\f{B^2}{s}\Big)dsd\tau + \int_{t_1}^{t_2} F(R,\tau)-\bar P + \f12 B^2d\tau\\
&\leqs \|\rho\|_\infty\int_0^T\|\g u\|_2^2dt + \int_0^T\|\g B\|_2^2dt+\Big|\int_{t_1}^{t_2} F(R,\tau)-\bar P+\f12 B^2d\tau\Big|\\
&\leqs \|\rho\|_\infty + 1 + (t_2-t_1).
\ea\]
where we had used the following fact
\be\label{fr-1-2}
\ba
\Big|\int_{t_1}^{t_2} F(R,\tau)d\tau\Big| \leqs \|\rho\|_\infty + 1 + (t_2-t_1).
\ea
\ee
 
To prove \eqref{fr-1-2}, we make use of Lemma \ref{FR}. Indeed, we had
\be\label{fr-1-1}
\ba
\int_{t_1}^{t_2} F(R,\tau)d\tau &=  \frac{1}{R^2} \Big[\int_{0}^{R} \rho u r^2 dr \big|_{\ts=t_2} + \int_{0}^{R} \rho u r^2 dr \big|_{\ts=t_1}\\
& + 2\int_{t_1}^{t_2}\int_{0}^{R} F(r, \ts) r dr d\ts - \int_{t_1}^{t_2}\int_{0}^{R} (\rho u^2 -B^2)r dr d\ts \Big]
\ea
\ee
Consequently,

\[\ba
\Big|\int_{t_1}^{t_2} F(R,\tau)d\tau\Big| &\leqs 1 + \Big|\int_{t_1}^{t_2}\int_{0}^{R} F(r, \ts) r dr d\ts\Big| +\| \rho\|_\infty\int_0^T\|\div u\|_2^2dt\\
&\leqs 1 + \| \rho\|_\infty + \Big|\int_{t_1}^{t_2}\int_{0}^{R} F(r, \ts) r dr d\ts\Big|.
\ea\]
Note that
\[\ba
\Big|\int_{0}^{R} F(r, \ts) r dr\Big| &= \Big|\int F(x, \ts) dx \Big|=\Big|\int \ld\div u - P + \bar P-\f12 B^2 dx\Big|\\
%&\lesssim \int \ld + \ld|\div u|^2 + P dx\Big|\\
&\lesssim \int \ld dx + \int \ld|\div u|^2dx + 1,
\ea\]
which gives that
\[
\Big|\int_{t_1}^{t_2}\int_{0}^{R} F(r, \ts) r dr d\ts \Big|
\leqs \int_{t_1}^{t_2}\int \ld dx + \int \ld|\div u|^2dx + 1d\ts
\leq 1+(t_2-t_1).
\]
with the help of the energy inequality and $\beta\le\gamma$.
This proves \eqref{fr-1-2}.

By Zlotnik's lemma,
\[
\theta + \xi \leqs \|\rho\|_\infty +1.
\]
\[
\rho^\beta \leqs \theta \leqs \|\xi\|_\infty + \|\rho\|_\infty +1 \leqs C(\a)R_T^{1+\a\beta} + R_T^{1+\gamma/2-\beta}.
\]
then,
\[
R_T^{\beta} \leqs C(\a)R_T^{1+\a\beta} + R_T^{1+\gamma/2-\beta}.
\]
Direct calculations show that
\be
R_T\le C
\ee
independent of time $T$ as long as 
\be
\max\{1,\frac{\gamma+2}{4}\}<\beta\le\gamma.
\ee
This finishes the proof of Proposition \ref{SZ1}.
\endproof

\section{Proofs of Theorem \ref{them1} and Theorem \ref{them2}}
With the key apriori estimates of upper bound of density at hand, the high order estimates and the proofs of Theorem \ref{them1} and Theorem \ref{them2} can be proved as a  standard procedure. See also \cite{chen2022,2016Huang} for more details.

\section*{Acknowledgements}
 X.-D. Huang is partially supported by NNSFC Grant Nos. 11971464, 11688101 and CAS Project for Young Scientists in Basic Research, Grant No.YSBR-031, National Key R\&D Program of China, Grant No.2021YFA1000800. W. YAN is partially supported by NNSFC Grant Nos. 11371064, 11871113. R. Yu is partially supported by Guangdong Provincial Natural Science Foundation, Grant No. 2020A1515110942.


\addcontentsline{toc}{section}{\refname}
\begin{thebibliography}{99}

\bibitem{BKM} J. Beal, T. Kato, A. Majda, Remarks on the breakdown of smooth solutions for the 3-D Euler equations, \textit{Comm. Math. Phys.} , {94}, 61-66 (1984).


\bibitem{BD1}
D. Bresch and B. Desjardins, Existence of global weak solutions for
a 2D viscous shallow water equations and convergence to the
quasi-geostrophic model, {\it Comm. Math. Phys.}, { 238}(1-2), 211-223 (2003).


\bibitem{bresch2007on}
D.~Bresch, B.~Desjardins, and D.~G\'{e}rard-Varet.
\newblock On compressible {N}avier-{S}tokes equations with density dependent
viscosities in bounded domains.
\newblock {\em J. Math. Pures Appl.}, 87(2):227--235, 2007.

\bibitem{chen2022} Y.Chen, B.Huang, X.Shi, 
\newblock Global strong and weak solutions to the
initial-boundary-value problem of two-dimensional
compressible mhd system with large initial data
and vacuum
\newblock {\em SIAM J. Math. Anal.} 54(3), 3817--3847, 2022.

\bibitem{2009Fan}
\newblock J.S. Fan, W.H. Yu.
\newblock Strong solution to the compressible magnetohydrodynamic equations with vacuum, \newblock {\rm Nonlinear Anal.}, 10:392--409,2009.

\bibitem{2022Li}
\newblock X.Y. Fan, J.X. Li and J. Li.
\newblock Global existence of strong and weak solutions to 2D compressible Navier-Stokes system in bounded domains with large data and vacuum.
\newblock {\em Arch.Rational Mech. Anal.}, 245:239--278, 2022.



\bibitem{2016Huang} X. D. Huang, J. Li, Existence and blowup behavior of global strong solutions to
the two-dimensional barotropic compressible Navier-Stokes system with vacuum
and large initial data, {\it J. Math. Pures Appl.}, 106(1), 123-154 (2016).

\bibitem{2022Huang} X. D. Huang, J. Li, Global Well-Posedness of Classical Solutions to the Cauchy Problem of Two-Dimensional Barotropic Compressible Navier--Stokes System with Vacuum and Large Initial Data,
{\it SIAM J. Math. Anal.} 54(3), 3192-3214, 2022.

\bibitem{huang2012global}
X. D.~Huang, J.~Li, and Z.~Xin.
\newblock Global well-posedness of classical solutions with large oscillations
and vacuum to the three-dimensional isentropic compressible {N}avier-{S}tokes
equations.
\newblock {\em Comm. Pure. Appl. Math.}, 65(4):549--585, 2012.


\bibitem{2023Huang}
\newblock X. D.~Huang, M.L.~Su, W. Yan and R.F.~Yu.
\newblock Global large strong solutions to the radially symmetric compressible Navier-Stokes
equations in 2D solid balls.
\newblock doi:https://arxiv.org/abs/2310.05040.

\bibitem{2023HXY} \newblock X. D.~Huang, W. Yan and Z. Xin.
\newblock On formation of finite time singularity of the two dimensional compressible magneto-hydrodynamic flows.
\newblock Preprint, 2023.

\bibitem{JWX} 
Q. Jiu, Y. Wang and Z. Xin.
\newblock Global classical solution to two-dimensional compressible Navier-Stokes equations with large data in $R^2$.
\newblock {\em Physica D},376/377: 180--194,2018.

\bibitem{1983Ka} 
\newblock S. Kawashima. 
\newblock Systems of a hyperbolic-parabolic composite type, with applications to the equations of magnetohydrodynamics,
\newblock {\rm PhD thesis}, Kyoto University, 1983.



\bibitem{2013lihl}
\newblock H.L. Li, X.Y. Xu, J.W. Zhang
\newblock  Global classical solutions to 3D compressible magnetohydrodynamic equations with large oscillations and vacuum
\newblock {\rm SIAM J. Math. Anal.}, 45:1356--1387,2013.


\bibitem{Lhlz2016}
H.-L. Li and X.W. Zhang.
\newblock Global strong solutions to radial symmetric compressible Navier-Stokes equations with free boundary.
\newblock {\em  J. Differential Equations},261(11):6341--6367,2016.


\bibitem{Kazhikhov} V. A. Vaigant, A. V. Kazhikhov, On the existence of global solutions of two-dimensional Navier-Stokes equations of a compressible viscous fluid. (Russian) {\it Sibirsk. Mat. Zh.}, 36 (1995), no. 6, 1283--1316, ii; translation in {\it Siberian Math. J.}, { 36}, no. 6, 1108-1141 (1995).


\bibitem{2015Mei-1}
Y.~Mei.
\newblock
Global classical solutions to the 2D compressible MHD equations with large data and vacuum.
\newblock {\em J. Differential Equations}258(9), 3304--3359,2015.

\bibitem{2015Mei-2}
Y.~Mei.
\newblock Corrigendum to ``Global classical solutions to the 2D compressible MHD equations with large data and vacuum'' [MR3317636].
\newblock {\em J. Differential Equations},258(9): 3360--3362,2015.

\bibitem{2022Shi}
\newblock Y.Z. Chen, B. Huang, X.D. Shi.
\newblock Global strong and weak solutions to the initial-boundary-value problem of two-dimensional compressible MHD system with large initial data and vacuum.
\newblock {\em SIAM.J. Math.Anal.}, 54(3):3817--3847,2022.


\bibitem{1995Kazhikhov}
V.~A. Vaigant and A.~V. Kazhikhov.
\newblock On existence of global solutions to the two-dimensional
{N}avier-{S}tokes equations for a compressible viscous fluid.
\newblock {\em Siberian Math. J.}, 36(6):1108--1141, 1995.


\bibitem{1972Vol}
\newblock A.I. Vol’pert, S.I. Hudjaev.
\newblock The Cauchy problem for composite systems of nonlinear differential equations,
\newblock {\rm Math.USSR-Sb}, 16: 517--544,1972.


\bibitem{zlot} A.A. Zlotnik.
\newblock Uniform estimates and stabilization of symmetric solutions of a system of quasilinear equations, {\em J. Differ. Equ.} 36 (2000) 701-716.

\end{thebibliography}
\end{document}